\documentclass[12pt]{article}
\usepackage{amsmath}
\usepackage{amsfonts}

\newtheorem{theorem}{Theorem}
\newtheorem{definition}{Definition}
\newtheorem{example}{Example}
\newtheorem{lemma}{Lemma}
\newtheorem{corollary}{Corollary}
\newtheorem{proposition}{Proposition}
\newtheorem{remark}{Remark}

\begin{document}

\title{Weights of modular forms on $\mathrm{SO}^{+}(2,l)$ and congruences between Eisenstein series and cusp forms of half-integral weight on $\mathrm{SL}_{2}$}
\author{Richard Hill}
\maketitle

\begin{abstract}
\noindent
Let $E$ be a level 1, vector valued Eisenstein series of half-integral weight, 
normalized so that the coefficients are all in $\mathbb{Z}$.
We show that there is a level one vector valued cusp form $f$ with the 
 same weight as $E$ and with coefficients in $\mathbb{Z}$,
 which is congruent to $E$ modulo the constant 
 term of $E$.
\end{abstract}

\section{Introduction}

Recall that the Ramanujan $\tau$ function satisfies the congruence
$$
	\tau(n) \equiv \sigma_{11}(n) \bmod 691,\quad
	\sigma_{11}(n)=\sum_{d|n}d^{11}.
$$
In another notation, this says that the weight 12 level 1 holomorphic Eisenstein series
$G_{12}$, is congruent modulo 691, to the cusp form $\Delta$.
In fact it is known, that for any holomorphic Eisenstein series $G$ of integral weight $k$ and arbitrary level $N$, and normalized so as to have integer coefficients,
 there is a cusp form $f$ of the same weight, with integer coefficients, such that
$$
	f\equiv G \bmod d,
$$
where $d$ denotes the constant term of $G$.
Congruences of this kind have applications in Bloch-Kato type conjectures,
 since the Galois representations corresponding to the cusp forms
 can be used to construct non-trivial elements in certain Selmer groups
 (this phenomenon is explained in the introduction to \cite{brown}).

In this paper we shall prove a similar congruence result (Theorem \ref{main})
 in the context of vector-valued forms of integral or half-integral weight on $\mathrm{Mp}_{2}(\mathbb{Z})$.
Only a few such congruences have been proved for half-integral weight
 (see \cite{AntoniadisKohnen}, \cite{koblitz}).
The method previously used in the half-integral weight cases was to
 deduce congruences in half-integral weight
 from those of integral weight using the Shimura correspondence
The proof in this paper is rather simple, given the prerequisites,
 and is entirely different from previous proofs of similar results;
it works equally well in both integral and half-integral weights.

To see how the congruence arises,
recall that a Borcherds lift is \`a-priori a rational weight meromorphic form on an orthogonal group $\mathrm{SO}^{+}(2,l)$ with respect to some arithmetic subgroup $\Gamma$
(see \cite{bruinier}).
Using an idea of Deligne \cite{deligne1} we show that if the congruence kernel of  the pre-image of $\Gamma$ in $\mathrm{Spin}(2,l)$ is trivial, then there only exist half-integral weight multiplier systems on $\Gamma$.
Hence, using Kneser's calculation \cite{kneser} of the congruence kernels of Spinor groups, we show (see Theorem \ref{orthogonalweights} below)
 that all forms on $\mathrm{SO}^{+}(2,l)$ have at most half-integral weight.
Consequently Borcherds lifts in fact have half-integral weight.
However, Borcherds lifts are constructed from certain weakly holomorphic vector valued forms on $\mathrm{SL}_{2}$.
The weight of the Borcherds lift is precisely half the constant term of the weakly holomorphic form.
As the weight is a half-integer, we deduce that the constant term is an integer.
This fact can be reinterpreted in terms of holomorphic forms on $\mathrm{SL}_{2}$ to give 
 the required congruence.

This paper was written while the author was a visiting Ulam professor in the University of Colorado at Boulder.
I'd like that thank the Department of Pure Mathematics at Boulder, especially Prof. Lynne Walling, for their hospitality.
I'd also like to thank Prof. Kevin Buzzard for a useful discussion.

\section{Statement of the result}

Let $\mathcal{H}$ denote the upper half-plane with the usual action of $\mathrm{SL}_{2}(\mathbb{R})$.
Following \cite{bruinier} we realise the metaplectic group $\mathrm{Mp}_{2}(\mathbb{R})$ as
 the set of pairs $(g,\phi)$,
 where $g=\left(\begin{matrix}a&b\\c&d\end{matrix}\right)\in \mathrm{SL}_{2}(\mathbb{R})$
 and $\phi:\mathcal{H}\to\mathbb{C}^{\times}$ is a continuous function
such that $\phi(\tau)^{2}=c\tau+d$. The group operation in $\mathrm{Mp}_{2}(\mathbb{R})$ is given by
$$
	(g,\phi)(h,\psi) := (gh, (\phi\circ h)\cdot \psi ).
$$
There is an obvious projection $\mathrm{Mp}_{2}(\mathbb{R})\to \mathrm{SL}_{2}(\mathbb{R})$
and there is a short exact sequence
$$
	1 \to \mu_{2} \to \mathrm{Mp}_{2}(\mathbb{R}) \to \mathrm{SL}_{2}(\mathbb{R}) \to 1,
$$
where $\mu_{2}=\{1,-1\}$ is embedded in $\mathrm{Mp}_{2}(\mathbb{R})$ by
 $\epsilon\mapsto (I_{2},\epsilon)$.
The centre of $\mathrm{Mp}_{2}(\mathbb{R})$ is $\mu_{2}$.

We shall write $\mathrm{Mp}_{2}(\mathbb{Z})$ for the preimage of $\mathrm{SL}_{2}(\mathbb{Z})$ in $\mathrm{Mp}_{2}(\mathbb{R})$.
This group $\mathrm{Mp}_{2}(\mathbb{Z})$ is generated by the following two elements
$$
	T
	=
	\left(
	\left(\begin{matrix}1&1\\0&1\end{matrix}\right),1
	\right),\quad
	S
	=
	\left(
	\left(\begin{matrix}0&1\\-1&0\end{matrix}\right),\sqrt{\tau}
	\right).
$$

Let $L$ be lattice over $\mathbb{Z}$ with a $\mathbb{Z}$-valued quadratic form $q:L\to\mathbb{Z}$.
We shall assume throughout that $q$ has signature $(2,l)$ with $l\ge 3$.
We shall also write $(\cdot,\cdot)$ for the corresponding symmetric bilinear form
and $L'$ for the dual lattice:
$$
	q(v)=\frac{1}{2}(v,v),\qquad
	L'=\{v\in L\otimes \mathbb{Q} : (v,L)\subseteq\mathbb{Z}\}.
$$
We shall consider the vector space $\mathbb{C}[L'/L]$ of complex-valued functions
on $L'/L$.
For $\alpha\in L'/L$ we let $[\alpha]$ denote the function
 which takes value $1$ at $\alpha $ and is zero elsewhere.
Thus the functions $[\alpha]$ form a basis of $\mathbb{C}[L'/L]$.
There is a representation $\varrho$ of $\mathrm{Mp}_{2}(\mathbb{Z})$
 on $\mathbb{C}[L'/L]$, which is a special case of the Weil representation (see \cite{weil}).
The representation $\varrho$ is discussed in more detail in
 \cite{bruinier} and \cite{shintani}).
On the generators $S$ and $T$, $\varrho$ is given by
$$
	\varrho(S)[\alpha]
	=
	\frac{1}{\mathfrak{g}}
	\sum_{\beta\in L'/L}
	\exp(-2\pi i (\beta,\alpha))
	\cdot[\beta]
	\qquad
	\varrho(T)[\alpha]
	=
	\exp(2\pi i q(\alpha))\cdot[\alpha],
$$
where $\mathfrak{g}$ denotes the following Gauss sum:
$$
	\mathfrak{g}
	=
	\sum_{\alpha\in L'/L}\exp(2\pi i q(\alpha)).
$$
Let $k\in\frac{1}{2}\mathbb{Z}$.
By a weight $k$, $\mathbb{C}[L'/L]$-valued modular form, we shall mean a holomorphic function
$f:\mathcal{H} \to \mathbb{C}[L'/L]$, such that
\begin{itemize}
	\item[(A)]
	for all $(\gamma,\phi)\in \mathrm{Mp}_{2}(\mathbb{Z})$ we have
	$$
	f(\gamma \tau)
	=
	\phi(\tau)^{2k}\cdot \varrho(\gamma,\phi) f(\tau),
	$$
	\item[(B)]
	$f$ is bounded on the region $\{\tau\in\mathbb{C} : \Im\tau >1\}$.
\end{itemize}
The first condition
implies that $f$ has a Fourier expansion of the form
$$
	f(\tau)
	=
	\sum_{\alpha\in L'/L}
	\sum_{n\in \mathbb{Z}+q(\alpha)}
	c(f,n, \alpha)q^{n}
	\cdot [\alpha],
	\qquad
	q=\exp(2\pi i \tau),\qquad
	c(f,n, \alpha)\in \mathbb{C},
$$
and the second condition implies that $c(f,n, \alpha)=0$ for $n<0$.
The form $f$ is called a cusp form if $c(f,n, \alpha)=0$ whenever $n=0$.
We shall write $\mathcal{M}_{k,L}$ for the complex vector space of weight $k$ forms,
 and $\mathcal{S}_{k,L}$ for the subspace of cusp forms.

\begin{remark}
	Suppose $f\in \mathcal{M}_{k,L}$, and expand $f$ in terms of the basis
	$[\alpha]$:
	$$
		f(\tau)
		=
		\sum_{\alpha\in L'/L}
		f_{\alpha}(\tau)\cdot [\alpha].
	$$
	Then each component $f_{\alpha}$ is a (scalar valued) holomorphic modular
	form of weight $k\in\frac{1}{2}\mathbb{Z}$ for the principal congruence subgroup
	$\Gamma(N)$, where $N\cdot L'\subseteq L$.
\end{remark}

We need the following result of McGraw \cite{mcgraw}:

\begin{theorem}
	\label{rationality}
	The space $\mathcal{M}_{k,L}$ has a basis of modular forms,
	all of whose Fourier expansions have only integer coefficients.
	Similarly, $\mathcal{S}_{k,L}$ has a basis of modular forms,
	all of whose Fourier expansions have only integer coefficients.
\end{theorem}

For example, let $k=\frac{1}{2}\mathrm{rank}(L)=1+\frac{l}{2}$ and define:
$$
	E(\tau)
	=
	\sum_{(\gamma,\phi)\in \tilde\Gamma_{\infty}\backslash \mathrm{Mp}_{2}(\mathbb{Z})}
	\phi(\tau)^{-2k}\varrho^{*}(\gamma,\phi)^{-1}[0]
$$
This Eisenstein series converges since $k>2$,
 and is an element of $\mathcal{M}_{1+l/2,L}$.
Bruinier and Kuss have shown that the Fourier coefficients $c(E,\beta,n)$
 of $E$ are in $\mathbb{Q}$ (see Corollary 8 of \cite{bruinierkuss}).
By Theorem \ref{rationality} (or by the above remark, together with results of
 \cite{serrestark})
 there is a non-zero integer $d$ such that $dE$ has integer coefficients.

The purpose of this paper is to prove the following:

\begin{theorem}
	\label{main}
	Let $E$ be the weight $k>2$ Eisenstein series defined above.
	Choose $d\in\mathbb{N}$ such that $d\cdot E$ has integer Fourier coefficients.
	Then there is a weight $k$ cusp for $f$
	 with integer Fourier coefficients $c(f,n,\alpha)$,
	 such that each coefficient satisfies the congruence:
	$$
		c(f,n, \alpha) \equiv c(d\cdot E,n, \alpha)\bmod d.
	$$
\end{theorem}

\begin{example}
	Let $q:\mathbb{Z}^{5}\to\mathbb{Z}$ be given by
	$$
		q\left(\begin{matrix}x_{1}\\\vdots\\x_{5}\end{matrix}\right)
		=
		x_{1}^{2}+x_{2}^{2}-x_{3}^{2}-x_{4}^{2}-x_{5}^{2}.
	$$
	In this case the space $\mathcal{S}_{5/2,L}$ of cusp forms is $\{0\}$.
	Therefore in this case the Eisenstein series
	 $E$ has coefficients in $\mathbb{Z}$.
\end{example}

\section{Weakly holomorphic modular forms}

By a \emph{weakly holomorphic} $\mathbb{C}[L'/L]$-valued modular form we shall mean
 (following the terminology of \cite{bruinierfunke})
a holomorphic function $f:\mathcal{H}\to \mathbb{C}[L'/L]$ satisfying condition (A) above,
and whose Fourier expansion has at most finitely many non-zero terms $c(f,\alpha,n)$ with $n<0$.
We shall write $\mathcal{M}_{k,L}^{!}$ for the space of weight $k$ weakly holomorphic forms.
Given an $f\in \mathcal{M}_{k,L}^{!}$, we define the \emph{principle part} of $f$ to be the
 negative part of its Fourier expansion:
$$
	\sum_{\alpha\in L'/L}\;\;
	\sum_{n\in \mathbb{Z}+q(\alpha),\; n<0}
	c(f, \alpha,n)\exp(2\pi i n \tau)\cdot [\alpha].
$$
It turns out that there exist weakly holomorphic forms with negative weight.
We shall be particularly interested in forms of weight $1-\frac{l}{2}$.

The following theorem of Borcherds
 (Theorem 3.1 of \cite{borcherds2} or Theorem 1.17 of \cite{bruinier})
 determines which Fourier polynomials arise
 as principle parts of weakly holomorphic forms:

\begin{theorem}
	\label{principalpart}
	Let $k\in \frac{1}{2}\mathbb{Z}$ and consider a Fourier polynomial
	of the form
	$$
		p(\tau)
		=
		\sum_{\alpha\in L'/L}\;\;
		\sum_{\begin{array}{c}
			n\in\mathbb{Z} +q(\alpha),\\
			-N\le n<0
		\end{array}}
		c(\alpha,n)\exp(2\pi i n \tau)\cdot [\alpha].
	$$
	There exists a weakly holomorphic modular form in $\mathcal{M}_{1-k,L}^{!}$
	 with principle part $p(\tau)$ if and only if for every
	 holomorphic cusp form $f\in \mathcal{S}_{k,L}$
	the following holds:
	$$
		\sum_{\alpha\in L'/L}\;\;
		\sum_{\begin{array}{c}
			n\in\mathbb{Z} +q(\alpha),\\
			-N\le n<0
		\end{array}}
		c(\alpha,n) \cdot c(f,\alpha,-n)
		=
		0.
	$$
\end{theorem}

If $f\in \mathcal{M}_{1-l/2,L}^{!}$ is a weakly holomorphic form,
 then its constant coefficient $c(0,0)$ may be expressed
 in terms of the principal part as follows
 (see Remark 3.23 of \cite{bruinier}, ignoring the misprint there):
\begin{equation}
	\label{constantterm}
	c(f,0,0)
	=
	-\frac{1}{2}
	\sum_{\alpha\in L'/L}\;\;
	\sum_{\begin{array}{c}
		n\in\mathbb{Z} +q(\alpha),\\
		n>0
	\end{array}}
	c(E,\alpha,n) \cdot c(f,\alpha,-n),
\end{equation}
where $c(E,\alpha,n)$ are the coefficients of the Eisenstein series $E$.

\section{Borcherds lifts to the Orthogonal group}

Recall that we have a lattice $L$ with a quadratic form $q$
 with signature $(2,l)$ ($l\ge 3$).
For a field $F$, we extend the bilinear form $(\cdot,\cdot)$ and the quadratic form $q$
 from $L$ to quadratic and $F$-bilinear form on $L\otimes_{\mathbb{Z}} F$.
Define
$$
	\mathrm{SO}(2,l)
	=
	\{\sigma\in \mathrm{SL}(V\otimes\mathbb{R}) : q\circ \sigma =q\}.
$$
As a Lie group, $\mathrm{SO}(2,l)$ has 2 connected components.
The connected component of the identity (which is the kernel
 of the spinor norm, see \cite{hahnomeara})
 will be denoted $\mathrm{SO}^{+}(2,l)$.
 
To make thinks more precise we shall choose a basis $\{b_{1},\ldots,b_{2+l}\}$
for $L\otimes\mathbb{R}$ such that
$$
	(b_{i},b_{j})
	=
	\begin{cases}
		0&\hbox{if $i\ne j$,}\\
		1&\hbox{if $i= j\le 2$,}\\
		-1&\hbox{if $i= j> 2$.}
	\end{cases}
$$
With this choice of basis, we shall regard $\mathrm{SO}(2)\oplus\mathrm{SO}(l)$ as a subgroup
 of $\mathrm{SO}(2,l)$ by the embedding
$$
	(A,B)
	\mapsto
	\left(\begin{matrix}
		A & 0\\
		0 & B
	\end{matrix}\right).
$$

For $Z\in L\otimes\mathbb{C}$ we let $\bar Z$ denote the complex conjugate
 with complex conjugation acting on $\mathbb{C}$.
Define
$$
	\tilde{\mathcal{K}}
	=
	\{
	Z\in L\otimes \mathbb{C} :
	(Z,Z)=0 \hbox{ and } (Z,\bar Z)>0
	\},
	\quad
	\mathcal{K}
	=
	\mathbb{P}(\mathcal{K}),
$$
where $\mathbb{P}(\mathcal{K})$ denotes the image of $\mathcal{K}$ in
projective space.

It is clear the the standard actions of $\mathrm{SO}(2,l)$ on
 $\mathbb{C}^{2+l}$ and $\mathbb{P}^{1+l}(\mathbb{C})$ preserve the subsets $\tilde{\mathcal{K}}$ and $\mathcal{K}$.
Furthermore one can show that the action of $\mathrm{SO}(2,l)$ on $\mathcal{K}$ is transitive.
Also, if $z\in L\otimes\mathbb{R}$ is a non-zero isotropic vector then one can check that for all
 $Z\in\tilde{\mathcal{K}}$,
$$
	(Z,z)\ne 0.
$$
Consider the vector
$$
	Z_{0}
	=
	b_{1}+ib_{2}
	=
	\left(\begin{matrix} 1\\i\\0 \\\vdots\\0\end{matrix}\right)\in \tilde{\mathcal{K}}.
$$
One easily checks that the stabilizer of $[Z_{0}]$ in $\mathrm{SO}(2,l)$ is
$\mathrm{SO}(2)\oplus\mathrm{SO}(l)$.
Therefore $\mathcal{K}$ is homeomorphic to $\mathrm{SO}(2,l)/(\mathrm{SO}(2)\oplus\mathrm{SO}(l))$.
Furthermore if we define $\mathcal{K}^{+}$ to be the connected component of $[Z_{0}]$
 in $\mathcal{K}$ then we see that $\mathcal{K}^{+}$ is homeomorphic to
 $\mathrm{SO}^{+}(2,l)/(\mathrm{SO}(2)\oplus\mathrm{SO}(l))$.
As $\mathrm{SO}(2)\oplus\mathrm{SO}(l)$ is a maximal compact subgroup of $\mathrm{SO}^{+}(2,l)$,
 it follows that $\mathcal{K}^{+}$ is a symmetric space, and hence contractible.

We now fix an isotropic vector $z\in L$.

\begin{definition}
	For $\sigma\in\mathrm{SO}^{+}(2,l)$ and $Z\in \tilde{\mathcal{K}}^{+}$ we define
	$$
		j(\sigma,Z)
		=
		\frac{(\sigma(Z),z)}{(Z,z)}.
	$$
	This is non-zero, and clearly depends only on the image of $Z$ in $\mathcal{K}$.
\end{definition}

It follows immediately that $j$ is a factor of automorphy, i.e. it satisfies the relation:
$$
	j(\sigma_{1}\sigma_{2},Z)
	=
	j(\sigma_{1},\sigma_{2}Z)
	j(\sigma_{2},Z).
$$
Let $r=p/q\in\mathbb{Q}$.
Since the set $\mathcal{K}^{+}$ is contractible,
 we may choose, for each $\sigma\in\mathrm{SO}^{+}(2,l)$,
 a branch $j(\sigma,Z)^{r}$ of the $r$-th power of $j(\sigma,Z)$,
 which is continuous in $Z$.
We shall choose these branches in such a way that $j(\sigma,Z)^{r}$ is a Borel 
 measurable function of $(\sigma,Z)$ (we could, for example, choose it to be continuous
 on each double coset in the Bruhat decomposition of $\mathrm{SO}^{+}(2,l)$).

\begin{definition}
	\label{wrdefn}
	For $\sigma_{1},\sigma_{2}\in \mathrm{SO}^{+}(2,l)$
	we define
	$$
		w_{r}(\sigma_{1},\sigma_{2})
		=
		\frac{j(\sigma_{1}\sigma_{2},Z)^{r}}
		{j(\sigma_{1},\sigma_{2}Z)^{r}j(\sigma_{2},Z)^{r}}.
	$$
	This is a $q$-th root of unity and
	 depends only on $\sigma_{1}$, $\sigma_{2}$, and not on $Z$.
\end{definition}

The function $w_{r}$ is a Borel measurable function on
 $\mathrm{SO}^{+}(2,l)\times \mathrm{SO}^{+}(2,l)$
 and satisfies the 2-cocycle relation:
$$
	w_{r}(\sigma_{1},\sigma_{2})
	w_{r}(\sigma_{1}\sigma_{2},\sigma_{3})
	=
	w_{r}(\sigma_{1},\sigma_{2}\sigma_{3})
	w_{r}(\sigma_{2},\sigma_{3}).
$$
Up to a measurable coboundary, $w_{r}$ is independent of the choice of the branches
 $j(\sigma,Z)^{r}$.
It therefore determines a central extension of $\mathrm{SO}^{+}(2,l)$ by $\mu_{q}$.
We shall now describe this covering group in another way.

Now regard $\mathrm{SO}(2)\oplus\mathrm{SO}(l)$ as a subgroup of $\mathrm{SO}^{+}(2,l)$
 in the obvious way.
We shall calculate the restriction of $w_{r}$ to these two subgroups.
Note that since $\mathrm{SO}(2)$ is isomorphic to $\mathbb{R}/\mathbb{Z}$ it has a unique connected $q$-fold
 central extension:
$$
	\begin{array}{ccccccccc}
		1& \to& \mu_{q}& \to& \mathrm{SO}(2)& \to& \mathrm{SO}(2)& \to& 1\medskip\\
		&&&&\sigma&\mapsto & \sigma^{q}
	\end{array}
$$

\begin{lemma}
	\label{covercalc}
	Let $r=p/q$ with $p$ and $q$ coprime.
	\begin{itemize}
		\item
		The restriction of $w_{r}$ to $\mathrm{SO}(2)$
		 corresponds to a connected $n$-fold covering group of $\mathrm{SO}(2)$.
		\item
		The restriction of $w_{r}$ to $\mathrm{SO}(l)$ is a coboundary.
	\end{itemize}
\end{lemma}

\textit{Proof.} 
Recall that the vector $Z_{0}\in \tilde{\mathcal{K}}^{+}$
 is an eigenvector of every element of $\mathrm{SO}(2)\oplus\mathrm{SO}(l)$.
For $\sigma\in\mathrm{SO}(l)$ we have $\sigma Z_{0}=Z_{0}$.
Hence for $\sigma\in\mathrm{SO}(l)$ we have $j(\sigma,Z_{0})=1$,
 so we may choose $j(\sigma,Z_{0})^{r}=1$.
With this choice, $w_{r}(\sigma_{1},\sigma_{2})=1$ for all
 $\sigma_{1},\sigma_{2}\in \mathrm{SO}(l)$.

On the other hand for
$$
	\sigma_{\theta}
	=
	\left(\begin{matrix}
		\cos\theta&\sin\theta\\
		-\sin\theta&\cos\theta
	\end{matrix}\right)
	\in\mathrm{SO}(2),
$$
we have
$$
	\sigma_{\theta} Z_{0}
	=
	\left(\begin{matrix}
		\exp(i\theta)\\
		i\exp(i\theta)\\
		0\\\vdots\\0
	\end{matrix}\right)
	=
	\exp(i\theta) Z_{0}.
$$
Hence $j(\sigma_{\theta},Z_{0})=\exp(i\theta)$.
The result follows easily from this.
\hfill$\square$
\medskip

\begin{definition}
	Let $r\in \mathbb{Q}$ and let $\Gamma$ be an arithmetic subgroup of $\mathrm{SO}^{+}(2,l)$.
	By a weight $r$ multiplier system on $\Gamma$ we shall mean a function
	$$
		\chi:\Gamma\to \mathbb{C}^{\times},
	$$
	such that for all $\gamma_{1},\gamma_{2}\in\Gamma$,
	$$
		\chi(\gamma_{1}\gamma_{2})
		=
		w_{r}(\gamma_{1},\gamma_{2})
		\chi(\gamma_{1})\chi(\gamma_{2}).
	$$
\end{definition}

The homogeneous space $\mathcal{K}^{+}$ also has a complex structure
 (see \cite{bruinier}) which gives sense to the following definition and theorem.

\begin{definition}
	Let $\Gamma$ be an arithmetic subgroup of $\mathrm{SO}^{+}(2,l)$.
	By a weight $r$ meromorphic modular form with respect to $\Gamma$
	we shall mean a meromorphic function $\Psi:\mathcal{K}^{+}\to\mathbb{C}\cup\{\infty\}$,
	for which there is a weight $r$ multiplier system $\chi$ on $\Gamma$,
	satisfying for all $\gamma\in\Gamma$, $Z\in\mathcal{K}^{+}$,
	$$
		\Psi(\gamma Z)
		=
		\chi(\gamma)j(\gamma,Z)^{r}\Psi(Z).
	$$
\end{definition}

\begin{theorem}[Borcherds \cite{borcherds}]
	Let $f\in \mathcal{M}_{1-l/2,L}^{!}$ be a weakly holomorphic form,
	 whose Fourier coefficients $c(f,\delta,n)$ are integral for all $n<0$.
	Then there is a corresponding non-zero meromorphic form $\Psi_{f}$ on
	$\Gamma_{L}$ of weight $c(f,0,0)/2$.
	Here $\Gamma_{L}$ denotes the subgroup of elements in $\sigma\in\mathrm{SO}^{+}(2,l)$
	such that $\sigma L=L$. 
\end{theorem}

(This is the only part of Theorem 3.22 of \cite{bruinier} which we require).

\section{Covering groups of $\mathrm{SO}^{+}(2,l)$}

In this section, we shall show that pulling a modular form from $\mathrm{SO}^{+}(2,l)$ back to $\mathrm{Spin}(2,l)$ does not change its weight.

\begin{proposition}
	The fundamental group of $\mathrm{SO}^{+}(2,l)$ is isomorphic to $\mathbb{Z}\oplus(\mathbb{Z}/2)$.
	More precisely the obvious map
	$$
		\mathrm{SO}(2)\oplus \mathrm{SO}(l)\to \mathrm{SO}^{+}(2,l)
	$$
	induces an isomorphism
	$$
		\pi_{1}(\mathrm{SO}(2))\oplus \pi_{1}(\mathrm{SO}(l)) \to \pi_{1}(\mathrm{SO}^{+}(2,l))
	$$
\end{proposition}

\textit{Proof.} 
It follows from the Iwasawa decomposition of $\mathrm{SO}^{+}(2,l)$
is homotopic to its maximal compact subgroup $\mathrm{SO}(2)\oplus \mathrm{SO}(l)$.
\hfill$\square$
\medskip

The proposition implies that the universal cover
$\widetilde{\mathrm{SO}}^{+}(2,l)^{univ}$ fits into an exact sequence:
$$
	1 \to \mathbb{Z}\oplus (\mathbb{Z}/2) \to \widetilde{\mathrm{SO}}^{+}(2,l)^{univ} \to \mathrm{SO}^{+}(2,l)\to 1.
$$
Recall that any other connected covering group
 is of the form $\widetilde{\mathrm{SO}}^{+}(2,l)^{univ}/H$,
 for a unique subgroup $H\subset \mathbb{Z}\oplus(\mathbb{Z}/2)$.
The fundamental group of the covering group $\widetilde{\mathrm{SO}}^{+}(2,l)^{univ}/H$
 may be identified in an obvious way with $H$.

\begin{proposition}
	The double cover $\mathrm{Spin}(2,l)$ of $\mathrm{SO}^{+}(2,l)$
	corresponds to the subgroup of $\mathbb{Z}\oplus(\mathbb{Z}/2)$ generated by the element $(1,1)$.
\end{proposition}

\textit{Proof.} 
The subgroup corresponding to $\mathrm{Spin}(2,l)$ has index 2 in $\mathbb{Z}\oplus(\mathbb{Z}/2)$.
There are three such subgroups, which are generated as follows:
$$
	\langle (2,0),(0,1)\rangle,\quad
	\langle (1,0)\rangle,\quad
	\langle (1,1)\rangle.
$$
However the restrictions of $\mathrm{Spin}(2,l)$ to $\mathrm{SO}(2)$ and $\mathrm{SO}(l)$ are the non-trivial covers $\mathrm{Spin}(2)$ and $\mathrm{Spin}(l)$ respectively.
From this we deduce that neither $\pi_{1}(SO(2))$ nor $\pi_{1}(\mathrm{SO}(l))$ is contained in the subgroup.
Hence the only possibility is $\langle (1,1)\rangle$.
\hfill$\square$
\medskip

We shall identify $\pi_{1}(\mathrm{Spin}(2,l))$ with the subgroup $\langle (1,1)\rangle$,
 which is isomorphic to $\mathbb{Z}$.

For $q\in\mathbb{N}$ we let $\widetilde{\mathrm{SO}}^{+}(2,l)_{(q)}$ denote the $q$-fold cyclic cover of
 $\mathrm{SO}^{+}(2,l)$ corresponding to the subgroup
$$
	\langle (q,0),(0,1)\rangle.
$$
Our interest in this covering group is because Borcherds lifts are forms on this group.
We also let $\widetilde{\mathrm{Spin}}(2,l)_{(q)}$ denote the pullback of $\widetilde{\mathrm{SO}}^{+}(2,l)_{(q)}$ to $\mathrm{Spin}(2,l)$:
$$
	\begin{array}{ccc}
		\widetilde{\mathrm{Spin}}(2,l)_{(q)} & \to & \mathrm{Spin}(2,l)\medskip\\
		\downarrow&& \downarrow \medskip\\
		\widetilde{\mathrm{SO}}^{+}(2,l)_{(q)} & \to & \mathrm{SO}^{+}(2,l).
	\end{array}
$$

\begin{proposition}
	$\widetilde{\mathrm{Spin}}(2,l)_{(q)}$ is the unique connected $q$-fold cyclic cover
	of $\mathrm{Spin}(2,l)$.
\end{proposition}

\textit{Proof.} 
The cover $\widetilde{\mathrm{Spin}}(2,l)_{(q)}$ of $\mathrm{Spin}(2,l)$
corresponds to the intersection
 $\pi_{1}(\mathrm{Spin}(2,l))\cap \pi_{1}(\widetilde{\mathrm{SO}}^{+}(2,l)_{(q)})$.
This intersection is $\langle (q,q) \rangle$,
 which is the subgroup of $\pi_{1}(\mathrm{Spin}(2,l))$ of index $q$.
\hfill$\square$
\medskip

\begin{proposition}
	Let $r=p/q\in\mathbb{Q}$ with $p$ and $q$ coprime.
	The covering group of $\mathrm{SO}^{+}(2,l)$
	corresponding to the cocycle $w_{r}$
	(see Definition \ref{wrdefn})
	 is isomorphic to the connected $q$-fold cover
	 $\widetilde{\mathrm{SO}}^{+}(2,l)_{(q)}$.
\end{proposition}

\textit{Proof.} 
Let $\tilde G$ denote the covering group corresponding to $w_{r}$.
There is a canonical homomorphism $\widetilde{\mathrm{SO}}^{+}(2,l)^{univ}\to \tilde G$,
and by restriction a map $\varphi:\pi_{1}(\widetilde{\mathrm{SO}}(2,l)^{univ})\to \mu_{q}$.
$$
	\begin{array}{ccccccccc}
		1&\to&
		\mathbb{Z}\oplus(\mathbb{Z}/2)& \to&
		\widetilde{\mathrm{SO}}^{+}(2,l)^{univ} & \to&
		\mathrm{SO}^{+}(2,l)&\to& 1 \medskip \\
		&& \varphi\downarrow\;\;\; && \downarrow &&|| \medskip \\
		1&\to&
		\mu_{q}& \to&
		\tilde G & \to&
		\mathrm{SO}^{+}(2,l)&\to& 1. \\
	\end{array}
$$
We must show that the kernel of $\varphi$
 is precisely the subgroup $\langle (q,0),(0,1)\rangle$.
To calculate $\varphi$,
 it is enough calculate its restrictions to $\pi_{1}(\mathrm{SO}(2))=\mathbb{Z}$ and $\pi_{1}(\mathrm{SO}(l))=\mathbb{Z}/2$.

By Lemma \ref{covercalc}
we know that $\varphi$ is trivial on $\pi_{1}(\mathrm{SO}(l))$.
Furthermore by Lemma \ref{covercalc}, the preimage of $\mathrm{SO}(2)$ in $\tilde G$ is connected.
This implies that the induced map $\varphi:\pi_{1}(\mathrm{SO}(2))\to \mu_{q}$
 is surjective.
The result follows.
\hfill$\square$

\section{Fractional weight multiplier systems}

Let $G$ be a (real) connected Lie group with a connected cyclic cover
$$
 1 \to \mu_{n} \to \tilde G \to G \to 1.
$$
Here $\mu_{n}$ denotes the group of $n$-th roots of unity in $\mathbb{C}$.
Suppose we have an arithmetic subgroup $\Gamma\subset G$.
We shall discuss the following question:
\begin{itemize}
    \item[]
    does $\Gamma$ lift to a subgroup of $\tilde G$?
\end{itemize}
Let $\mathbb{C}^{1}$ denote the group of complex numbers with absolute value 
 1.
Suppose $w:G\times G \to \mu_{n}$ is a 2-cocycle representing
 the group extension $\tilde G$.
By a weight $w$ multiplier system on $\Gamma$,
 we shall mean a function $\chi: \Gamma \to \mathbb{C}^{1}$ such that
$$
 \chi(\gamma_{1}\gamma_{2})
 = 
 w(\gamma_{1},\gamma_{2}) \chi(\gamma_{1}) \chi(\gamma_{2}).
$$
In other words the image of $w$ in $Z^{2}(\Gamma, \mathbb{C}^{1})$
 is the coboundary $\partial\chi$.
If an arithmetic subgroup $\Gamma$ lifts to $\tilde G$
 then such a $\chi$ exists on $\Gamma$.
There is a partial converse to this:

\begin{proposition}
	\label{easy}
    If there is a weight $w$ multiplier system on
     an arithmetic subgroup $\Gamma\subset G$
     then there is an arithmetic subgroup $\Gamma_{0}\subset \Gamma$
     which lifts to $\tilde G$.
\end{proposition}

\begin{theorem}
	\label{weights}
	Let $G/\mathbb{R}$ be absolutely simple and simply connected
	and let $\tilde G\to G$ be a connected $n$-fold cyclic cover.
	Let $\Gamma$ be a congruence subgroup of $G$
	 such that every subgroup of finite index in $\Gamma$ is a congruence
	 subgroup.
	Furthermore, in the case that $G$ is a special unitary group,
	 assume that the construction of $\Gamma$ does not involve is a
	 non-abelian division algebra.
	If $\Gamma$ lifts to $\tilde G$ then $n\le 2$.
\end{theorem}	

Both these results are proved in \cite{hillkyoto}.
We shall apply the theorem in the case $G=\mathrm{Spin}(2,l)$.
To do this we need the congruence subgroup property for spinor groups:

\begin{theorem}[Kneser \cite{kneser}]
	\label{kneser}
	Let $L$ be a lattice over $\mathbb{Z}$ of signature $(2,l)$ with $l\ge 3$,
	and let $\Gamma$ be the preimage in $\mathrm{Spin}(2,l)$
	of the group of special orthogonal transformations of $L$.
	Then every subgroup of finite index in $\Gamma$ is a congruence subgroup.
\end{theorem}

\section{Proofs of the main results}

\begin{theorem}
	\label{orthogonalweights}
	Let $f$ be a non-zero weight $r\in\mathbb{Q}$ modular form
	on $\mathrm{SO}^{+}(2,l)$ with respect to some arithmetic subgroup
	of $\Gamma_{L}$.
	Then $r\in \frac{1}{2}\mathbb{Z}$.
\end{theorem}

\textit{Proof.} 
Let $r=p/q$ with $p$ and $q$ coprime.
Hence there is a weight $p/q$ multiplier system on $\Gamma_{L}$.
It follows by Proposition \ref{easy},
 that some arithmetic subgroup $\Gamma_{1}\subset\Gamma_{L}$
 lifts to the covering group $\widetilde{\mathrm{SO}}^{+}(2,l)_{(q)}$.
On the other hand, since $\mathrm{Spin}(2,l)$ is algebraic, it follows that
 some arithmetic subgroup $\Gamma_{2}\subset \Gamma_{L}$ lifts to $\mathrm{Spin}(2,l)$.
Hence the intersection $\Gamma=\Gamma_{1}\cap\Gamma_{2}$ lifts to
 $\widetilde{\mathrm{Spin}}(2,l)_{(q)}$.
 
Note now that the preimage of $\Gamma$ in $\mathrm{Spin}(2,l)$ is an arithmetic subgroup
 of $\mathrm{Spin}(2,l)$ which lifts to the connected $q$-fold cover $\widetilde{\mathrm{Spin}}(2,l)_{(q)}$.
Theorems \ref{weights} and \ref{kneser} now imply that $q\le 2$.
\hfill$\square$
\medskip

\begin{remark}
	We could not use Theorem \ref{weights} directly on $\mathrm{SO}^{+}(2,l)$,
	since this group is not algebraically simply connected, and its arithmetic subgroups
	are not all congruence subgroups.
\end{remark}

\begin{corollary}
	\label{integrality}
	Let $f\in \mathcal{M}_{1-l/2,L}^{!}$ be a weakly holomorphic form.
	If the Fourier coefficients $c(f,\alpha,n)$ are integral for all $n<0$.
	Then $c(f,0,0)\in\mathbb{Z}$.
\end{corollary}

\textit{Proof.} 
By Borcherds theorem together there is a non-zero form $\Psi_{f}$ of weight
$c(f,0,0)/2$ on the orthogonal group $\mathrm{SO}^{+}(2,l)$.
The result follows from the previous theorem.
\hfill$\square$
\medskip

In the case that $L$ is unimodular, we have $\mathbb{C}[L'/L]=\mathbb{C}$
and $\varrho$ is the trivial representation of $\mathrm{Mp}_{2}(\mathbb{Z})$.
In this case, $l$ is always even (in fact $2+l$ is a multiple of $4$),
 and our results refer to $\mathbb{C}$-valued, integral weight, level 1 forms.
We note the following generalization of Corollary \ref{integrality} in this context.

\begin{remark}
	Let $f$ be a weakly holomorphic, $\mathbb{C}$-valued modular form of level $1$ and weight $k\in 2\mathbb{Z}$
	with $k\le 0$.
	Let $f$ have the following Fourier expansion:
	$$
		f(\tau)
		=
		\sum_{n=-N}^{\infty} c(f,n)q^{n}.
	$$
	If $c(f,n)\in\mathbb{Z}$ for all $n<0$ then $c(f,n)\in\mathbb{Z}$ for all $n$.
\end{remark}

This is fairly easy to prove and seems to be well-known, although I don't know
where the proof has been written down.
To prove it one can use Corollary \ref{integrality}
 and proceed by induction on $n$, multiplying $f$ by $\Delta^{-1}$ at each step.

\begin{lemma}
	Fix $N\in\mathbb{N}$, and let $a(\alpha,n)$ be integers for
	$\alpha\in L'/L$ and $n=1,\ldots ,N$.
	If for every cusp form $f\in \mathcal{S}_{1+l/2,L}$ we have
	$$
		\sum a(\alpha,n)c(f, \alpha,n)=0
	$$
	then we must also have
	$$
		\sum a(\alpha,n)c(E, \alpha,n) \in \mathbb{Z},
	$$
	where $c(E, \alpha,n)$ are the coefficients of the Eisenstein series.
\end{lemma}

\textit{Proof.} 
Since $\sum a(\alpha,n)c(f, \alpha,n)=0$ for all cusp forms $f$,
it follows from Theorem \ref{principalpart} that there is
a weakly holomorphic form $g\in \mathcal{M}^{!}_{1-l/2,L}$ with principal part
$$
	c(g, \alpha,-n)=a(\alpha,n).
$$
By (\ref{constantterm}) the constant coefficient of $g$ is given by the formula
$$
	c(g,0,0)=\sum a(\alpha,n)c(E, \alpha,n).
$$
By the previous corollary this is an integer.
\hfill$\square$
\medskip

\begin{lemma}
	Let $N\in\mathbb{N}$.
	There is a cusp form $f\in \mathcal{S}_{1+l/2,L}$ with rational coefficients
	 such that for $0\le n\le N$ and $\alpha\in L'/L$,
	$$
		c(f, \alpha,n) \equiv c(E, \alpha,n) \bmod \mathbb{Z}.
	$$
\end{lemma}

\textit{Proof.} 
We can choose $N$ large enough so that
a modular form $f\in \mathcal{M}_{1+l/2,L}$ is determined by
the coefficients $c(f,\alpha,n)$ for $0\le n\le N$.
We shall write $\mathcal{M}_{1+l/2,L}(\mathbb{Q})$ and $\mathcal{S}_{1+l/2,L}(\mathbb{Q})$
 for the spaces of forms and cusp forms with rational coefficients.

Let $V$ be the group of functions $a:L'/L\times \{1,\ldots, N\}\to\mathbb{Z}$.
There is a (degenerate) pairing $V\times \mathcal{M}_{1+l/2,L}(\mathbb{Q})\to\mathbb{Q}$ given by
$$
	\langle a,f \rangle
	=
	\sum_{n=1}^{N} a(\alpha,n)c(f, \alpha,n).
$$
Let $V_{0}$ be the orthogonal complement of $\mathcal{M}_{1+l/2,L}(\mathbb{Q})$,
so we have a non-degenerate pairing $(V/V_{0})\times \mathcal{M}_{1+l/2,L}(\mathbb{Q})\to\mathbb{Q}$.
Let $W\subseteq V$ be the orthogonal complement of the space of cusp forms.
It follows from Theorem \ref{rationality} and our choice of $N$
 that $\mathcal{S}_{1+l/2,L}(\mathbb{Q})$ is the orthogonal complement of $W$.

The previous lemma shows that
$$
	\langle W, E \rangle
	\subseteq
	\mathbb{Z}.
$$
Since $W$ is a direct summand of $V$,
 it follows that there is a linear map $\lambda:V\to \mathbb{Z}$ such that
$$
	\lambda(w)=\langle w,E \rangle,
	\qquad
	w\in W.
$$ 
Since $\lambda$ is zero on $V_{0}$,
 there is a form $f\in \mathcal{M}_{1+l/2,L}(\mathbb{Q})$, such that
$$
	\langle v, f \rangle
	=
	\lambda(v),
	\qquad
	v\in V.
$$
Hence $E-f$ is in the orthogonal complement of $W$, so is a cusp form.
On the other hand, $\langle V,f \rangle=\lambda(V)\subset\mathbb{Z}$,
 so the coefficients $c(f,\alpha,n)$ are integers for $n\le N$.
\hfill$\square$
\medskip

\begin{theorem}
	There is a cusp form $f\in \mathcal{S}_{1+l/2,L}$ such that
	For all $n\in\mathbb{N}$ and $\alpha\in L'/L$,
	$$
		c(f, \alpha,n) \equiv c(E, \alpha,n) \bmod \mathbb{Z}.
	$$
\end{theorem}

\textit{Proof.} 
For the moment fix $N$.
We know that there is a cusp form $f$ such that for $n \le N$ and $\alpha\in L'/L$
we have $c(E, \alpha,n)-c(f, \alpha,n)\in \mathbb{Z}$.
Let $A_N$ be the set of all such $f$.
We clearly have $A_{1}\supseteq A_{2}\supseteq \ldots$.
Our aim is to show that the intersection of the sets $A_{N}$ is non-empty.

Clearly $A_N$ is a coset of a subgroup $B_N$ of $\mathcal{S}_{k,L}$.
In fact $B_N$ is the set of cusp forms whose first $N$ coefficients are integers.
As long as $N$ is large enough, a cusp form is determined by its first $N$ coefficients.
Hence $B_N$ is a $\mathbb{Z}$-module of finite rank.
As we have $B_{N+1} \subseteq B_N$, it follows that for large $r$ the rank of $B_r$ is stable.
Furthermore as the coefficients of cusp forms have bounded denominators
 (Theorem \ref{rationality} or \cite{serrestark}),
 it follows that for $N$ large enough we have $B_N=B_{N+1}= \ldots$.
Hence $A_N=A_{N+1}= \ldots$.
It follows that the intersection of all the sets $A_{N}$ is non-empty, which proves the theorem.
\hfill$\square$
\medskip

Theorem \ref{main} is equivalent to the above theorem.

\end{document}